# A compact implementation of a recently proposed strongly polynomial-time algorithm for the general LP problem

Samuel Awoniyi
Department of Industrial and Manufacturing Engineering
FAMU-FSU College of Engineering
2525 Pottsdamer Street
Tallahassee, FL 32310
E-mail: awoniyi@eng.famu.fsu.edu
ORCID# 0000 0001 7102 6257

## Abstract

This article presents a compact implementation of a recently proposed strongly polynomial-time algorithm for the general linear programming problem. Each iteration of the algorithm consists of applying a pair of complementary Gauss-Jordan (GJ) pivoting operations. In this compact implementation of the algorithm, the GJ pivoting operations are done inside a matrix that has half the size of the matrix used in the proposed algorithm. A numerical illustration is given.

## 1. Introduction

This article presents a relatively compact implementation of a recently proposed strongly polynomial-time algorithm for solving the general linear programming (LP) problem [6]. That algorithm utilizes basic LP duality theory to translate solving the general LP problem, having $k$ inequality constraints and $n$ variables, into solving a special system of equations in $R^{2(k+n)}$.

Each iteration of the proposed algorithm consists of two Gauss-Jordan (GJ) pivoting instances. The algorithm stops after at most 2(k+n) iterations. Comprehensive references are available at [1,2,3] on algorithms for the general LP problem.

This compact implementation of the proposed algorithm reduces each iteration to pivoting inside a (k+n+1)-by-(k+n+1) matrix in place of a (k+n+1)-by-2(k+n)+1 matrix. Each iteration of the compact implementation features a pair of complementary GJ pivoting transforming a skew-symmetric matrix in a manner that is quite instructive in its own right [5].

The rest of this article is organized as follows. Section 2 describes the compact implementation. A numerical illustration is given in Section 3. Remarks on related further work comprise Section 4.

.

## 2. Compact implementation description

We begin here by stating some notation and definitions that are needed for describing our compact implementation of the algorithm proposed in [6]. Thereafter, the steps of the compact implementation will be described.

### 2.1 Notation and definitions

As notation in this article, vectors are column vectors unless otherwise indicated. Vectors will be denoted by lower-case letters, and matrices by upper-case letters. Superscript $T$ will indicate vector or matrix transpose as usual, and $I_{(.)}$ is reserved for identity matrix of dimension indicated by (.).

We assume the general LP problem to be given in Neumann symmetric form, (P) below:

.

$$\left\{ \begin{array}{ll} \text{maximize} & f^T x \\ \text{subject to:} & Ax \le b, \\ & x \ge 0 \end{array} \right\} \cdots\cdots (P)$$

.

where $f$ is $n$-vector, $A$ is $k$-by-$n$ (numerical) matrix, $b$ is $k$-vector, and $x$ is $n$-vector of problem's variables.

From basic LP duality theory, solving (P) is equivalent to computing a $2(k+n)$-vector $z$ that solves the constrained system of linear equations, (Eq), stated below:

.

$$\left\{ \begin{array}{l} Mz = q, \\ z_j z_{(k+n+j)} = 0, \ \text{for } j = 1, \ldots, k+n \\ z \ge 0 \end{array} \right\} \cdots\cdots (Eq)$$

where

$$M = \begin{pmatrix} O & A & I_{(k)} & O \\ -A^T & O & O & I_{(n)} \\ -b^T & f^T & o^T & o^T \end{pmatrix} \text{and } q = \begin{pmatrix} b \\ -f \\ o \end{pmatrix}$$

.

As in the article [6], the (k+n+1)-by-(2k+2n+1) augmented matrix obtained by putting $M$ and $q$ together, with $q$ as column 2k+2n+1, will be denoted by $[M\ q]$.

Next, we define the concept of complementary GJ+ pivoting. It is an enhancement of complementary GJ pivoting. It will play a role in this compact implementation similar to the role played by complementary GJ pivoting in [6].

.

Definition of complementary GJ+ pivoting in column $l$ of given s-by-s matrix $S \equiv (S_{i,j})$ having $S_{l,l} \neq 0$ or $S_{s,l} \neq 0$.

A complementary GJ+ pivoting in column $l$ of given square matrix S is defined by the three matrix operations specified in (i), (ii), (iii) below:

(i) augment S by attaching the $l$-th unit vector as column s+1, thereby obtaining an s-by-(s+1) matrix, say Q;

(ii) next perform GJ pivoting with position $(l, l)$ in Q as pivot position, possibly after adding the last row, row s, to the $l$-th row;

(iii) next swap column $l$ and column s+1 in the transformed Q, and thereafter drop off resultant column s+1 (which is now the $l$-th unit vector) from transformed Q, thereby obtaining the desired transformation of S (having the same size as S).

## 2.2 Steps of this compact implementation

We begin here with an overview intended to aid some intuition. The overview ends with a flow chart. Thereafter, we will describe details of "initialization step", "general iteration step" and "stopping step".

### 2.2.1 Cursory overview of compact implementation

Each iteration of the algorithm proposed in [6] consists of applying a pair of complementary GJ pivoting inside the (k+n+1)-by-(2k+2n+1) matrix $[M\ q]$. In each iteration, the first GJ pivoting is called "MinorP" pivoting, and the second GJ pivoting is called "MajorP" pivoting.

The 2k+2n+1 columns of each $[M\ q]$ instance always include k+n unit vectors in $R^{k+n+1}$ that together form an identity matrix in $R^{k+n}$. By shedding those k+n unit vectors, we represent the information contained in each $[M\ q]$ instance with a (k+n+1)-by-(k+n+1) matrix, instead of a (k+n+1)-by-(2k+2n+1) matrix. Thus, each complementary GJ pivoting in $[M\ q]$ in the underlying algorithm (of [6]) corresponds to a GJ+ pivoting in a (k+n+1)-by-(k+n+1) matrix.

Accordingly, as a summary of this cursory overview of the iteration of both the proposed algorithm and its compact implementation, a flow chart connecting "Initialization", "MinorP instance", "MajorP instance" and "Stopping" is displayed here:

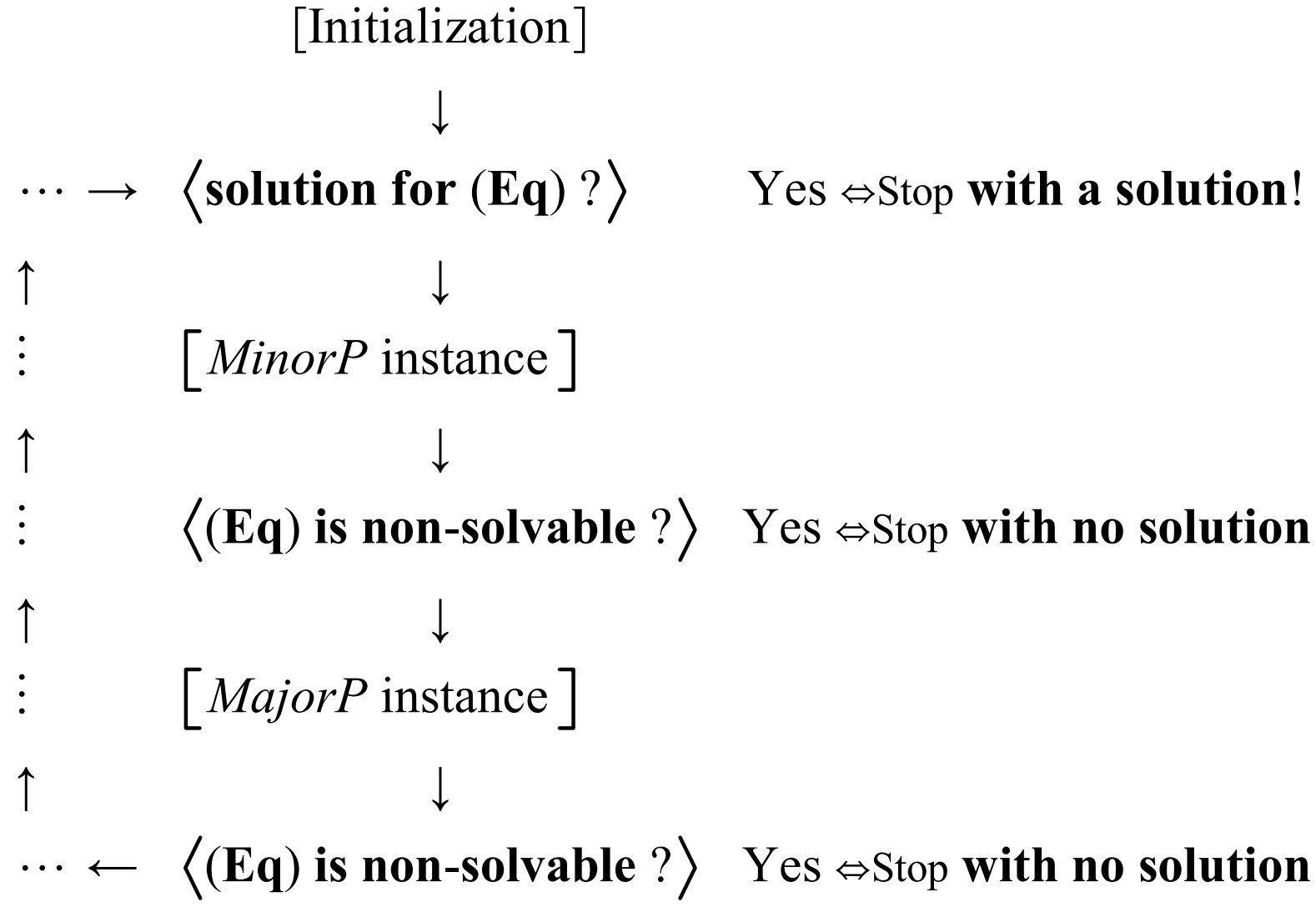

### 2.2.2 Initialization step for compact implementation

Define $P^{(0)}$ by

$$P^{(0)} \equiv \begin{pmatrix} \mathrm{O} & \mathrm{A} & \mathrm{b} \\ -\mathrm{A}^T & \mathrm{O} & -\mathrm{f} \\ -\mathrm{b}^T & \mathrm{f}^T & \mathrm{o} \end{pmatrix}$$

.

$$\equiv \begin{array}{|c|c|c|c|} \hline \mathrm{p}^0_{1,1} & \cdots & \mathrm{p}^0_{1,k+n} & \mathrm{p}^0_{1,k+n+1} \\ \hline \vdots & \ddots & \vdots & \vdots \\ \hline \mathrm{p}^0_{k+n,1} & \cdots & \mathrm{p}^0_{k+n,k+n} & \mathrm{p}^0_{k+n,k+n+1} \\ \hline \mathrm{p}^0_{k+n+1,1} & \cdots & \mathrm{p}^0_{k+n+1,k+n} & 0 \\ \hline \end{array}$$

.

If $(\mathrm{p}^0_{1,k+n+1}, \ldots, \mathrm{p}^0_{k+n,k+n+1}) \geq 0$, then the algorithm is terminated there, because (Eq) then has a trivial solution. Otherwise, this compact implementation next proceeds into the iterations, by setting iteration counter $i = 1$, and going to iteration 1.

The following table, which will be referred to as "Column Selection Record" (CSR), is initialized at this "Initialization step".

| | Z | P |
|---|---|---|
| 1 | | |

CSR will maintain a record of columns of $P^{(0)}$ that are 'nominated' by the iterations for inclusion in a basis matrix for a solution of (Eq) (in accordance with Lemma 6.1 in [6]). CSR will also be utilized in 'extracting at the end' a solution of LP problem (P) and its dual problem.

As a part of this initilization step, the indices of columns of $P^{(0)}$ that are known to be 'included in' a basis matrix for a solution of (Eq) are to be put in a set denoted as $\Pi$ (adopting a notation of [6]). The index set $\Pi$ will be updated as needed inside each iteration.

### 2.2.3 General iteration step for compact implementation

By virtue of the article [5], each iteration utilizes a pair of complementary GJ+ pivoting operations to produce two matrices, $Z^{(i)}$ & $P^{(i)}$ for the $i$-th iteration, $i = 1,2,\ldots.$ Each $Z^{(i)}$ corresponds to a MinorP pivoting (notation of [6]), and each $P^{(i)}$ corresponds to a MajorP pivoting.

To begin iteration 1 of this compact implementation, the matrix $Z^{(1)}$ is obtained by applying a GJ+ pivoting operation to $P^{(0)}$ (specified above under "Initialization"). Thereafter, the matrix $P^{(1)}$ is obtained by applying a GJ+ pivoting operation to $Z^{(1)}$.

In the $i$-th iteration, for $i = 2,3,\ldots,$ the matrix $Z^{(i)}$ will be obtained by applying a GJ+ pivoting operation to the matrix $P^{(i-1)}$, and the matrix $P^{(i)}$ in turn will be obtained by applying a GJ+ pivoting operation to the matrix $Z^{(i)}$.

More details of this $P^{(i-1)} \rightarrow Z^{(i)} \rightarrow P^{(i)}$ sequence are as follows. Towards obtaining $Z^{(i)}$ from $P^{(i-1)}$, for iteration $i = 2,3,\ldots,$ let

$$P^{(i-1)} \equiv \begin{array}{|c|c|c|c|} \hline \mathrm{p}^{i-1}_{1,1} & \cdots & \mathrm{p}^{i-1}_{1,k+n} & \mathrm{p}^{i-1}_{1,k+n+1} \\ \hline \vdots & \ddots & \vdots & \vdots \\ \hline \mathrm{p}^{i-1}_{k+n,1} & \cdots & \mathrm{p}^{i-1}_{k+n,k+n} & \mathrm{p}^{i-1}_{k+n,k+n+1} \\ \hline \mathrm{p}^{i-1}_{k+n+1,1} & \cdots & \mathrm{p}^{i-1}_{k+n+1,k+n} & 0 \\ \hline \end{array}$$

.

The last row and the last column of $P^{(i-1)}$ are related in a very special way explained in [5]. Therefore, by virtue of Lemma 1 in [5], we can assume, without loss of generality, that $\mathrm{p}^{i-1}_{k+n+1,j} > 0$ if $\mathrm{p}^{i-1}_{j,k+n+1} < 0,$ and $\mathrm{p}^{i-1}_{k+n+1,j} < 0$ if $\mathrm{p}^{i-1}_{j,k+n+1} > 0,$ for column/row index $j = 1,\ldots,k+n.$ The remaining details of how to obtain $Z^{(i)}$ from $P^{(i-1)}$ are given in the following 'Interpretation of MinorP instance" box.

*Interpretation of MinorP instance from [6]*

Define the set of column indices $C^{(\mathrm{i})}$ by

$$C^{(\mathrm{i})} = \{j \text{ such that } \mathrm{p}^{i-1}_{k+n+1,j} > 0\}, \text{ iteration } i = 1,2,\ldots$$

Define $L^{(\mathrm{i})}$ by $L^{(\mathrm{i})} \equiv C^{(\mathrm{i})} \backslash \Pi,\ i = 1,2,\ldots$

If $L^{(\mathrm{i})}$ contains only one element, say $\widehat{j}$, then put $\widehat{j}$ in $\Pi$, and perform complementary GJ+ pivoting in column $\widehat{j}$.

Otherwise order $L^{(\mathrm{i})}$ in ascending order of $\mathrm{p}^{i-1}_{k+n+1,j}$, and do:

(a) let $j^*$ be the first $j \in \mathrm{L}^{(\mathrm{i})}$ that's not entered in P column of current CSR. Perform complementary GJ+ pivoting in column $j^*$, and let the resultant matrix be the desired $Z^{(i)}$, and record $j^*$ under Z column in row $i$ of CSR. (It's advisable, but not necessary, to avoid repeating a $j$ in the Z column of CSR);

(b) but if such a $j^*$ (as specified in (a)) does not exist in $L^{(\mathrm{i})}$, then perform *separately* a pair of "complementary GJ+ pivoting in a column that is indexed in $L^{(\mathrm{i})}$, followed by a corresponding MajorP pivoting (see below)" until one such *separate* pair of GJ+ pivotings results in a solution of (Eq);

(c) but if such a successful column index (as specified in (b)) does not exist in $L^{(\mathrm{i})}$, then declare the conclusion that problem (Eq) has no solutions.

.

Towards obtaining $P^{(i)}$ from $Z^{(i)}$, for $i = 2,\ldots,$ let

.

$$Z^{(i)} \equiv \begin{array}{|c|c|c|c|}\hline z^i_{1,1} & \cdots & z^i_{1,k+n} & z^i_{1,k+n+1} \\ \hline \vdots & \ddots & \vdots & \vdots \\ \hline z^i_{k+n,1} & \cdots & z^i_{k+n,k+n} & z^i_{k+n,k+n+1} \\ \hline z^i_{k+n+1,1} & \cdots & z^i_{k+n+1,k+n} & z^i_{k+n+1,k+n+1} > 0 \\ \hline \end{array}$$

.

Regarding whether $z^i_{k+n+1,k+n+1} > 0$ in the last row of $Z^{(i)}$, it may be necessary to (implicitly) multiply the last row of $Z^{(i)}$ by -1, in order to ensure that $z^i_{k+n+1,k+n+1} > 0$. The remaining details of how to obtain $P^{(i)}$ from $Z^{(i)}$ are given in the following "Interpretation of MajorP instance" box.

.

*Interpretation of MajorP instance from [6]*

Define the set of column indices $C^{(\mathrm{i})}$ by

$C^{(\mathrm{i})} = \{\mathrm{j} \le \mathrm{k}+\mathrm{n}$ such that $\mathrm{z}^{i}_{k+n+1,j} > 0\}$, iteration $\mathrm{i} = 1,2,\ldots$

Define $\mathrm{L}^{(\mathrm{i})}$ by $\mathrm{L}^{(\mathrm{i})} \equiv C^{(\mathrm{i})} \backslash \Pi$, $i = 1,2,\ldots$

If $L^{(\mathrm{i})}$ contains only one element, say $\widehat{j}$, then put $\widehat{j}$ in $\Pi$ and perform the complementary GJ+ pivoting in column $\widehat{j}$; let the resultant matrix be the desired $P^{(i)}$, and record $\widehat{j}$ under P column in row $i$ of CSR.

Otherwise order $\mathrm{L}^{(\mathrm{i})}$ in descending order of $\mathrm{z}^{i}_{k+n+1,j}$, and do:

(a) let $j^*$ be next $j$ in $L^{(\mathrm{i})}$ that is not entered under P column in current CSR. Perform complementary GJ+ pivoting in column $j^*$, and let the resultant matrix be the desired $P^{(i)}$; record $j^*$ under the P column in row $i$ of CSR;

(b) but if such a $j^*$ (as described in (a)) does not exist in $L^{(\mathrm{i})}$, then perform a complementary GJ+ pivoting *separately* in every column that is indexed by $L^{(\mathrm{i})}$, until one such separate pivoting results in a solution of (Eq);

(c) but if such a successful column index (as described in (b)) does not exist in $L^{(\mathrm{i})}$, then declare the conclusion that problem (Eq) has no solutions.

### 2.2.4 Stopping step for compact implementation

There are two types of "stopping" – the case when a solution for (Eq) is found, and the case when there is clear evidence that (Eq) has no solutions.

*Case 1: A solution of (Eq) is found.* A solution of (Eq) is indicated in a $P^{(i)}$ instance by having (in the last column, column $k+n+1$)

$$(\mathrm{p}^{i}_{1,k+n+1},\ldots,\mathrm{p}^{i}_{k+n,k+n+1}) \ge 0 \text{ along with } \mathrm{p}^{i}_{k+n+1,k+n+1} = 0.$$

To obtain solutions of corresponding LP problem (P) and its dual problem, which one gets from the first $k+n$ components of a solution, say $y^*$, of (Eq) : for $j = 1,\ldots,k+n$, if $j$ has been entered into "Column Selection Record" (CSR) *an odd number of times*, then

$$\text{set } y^*_j = \mathrm{p}^{i}_{j,k+n+1}; \text{ otherwise set } y^*_j = 0.$$

*Case 2: There is clear evidence that (Eq) has no solutions.* A lack of solutions for (Eq) is indicated in a $Z^{(i)}$ by having (in the last row, row $k+n+1$)

$$(z^{i}_{1,k+n+1},\ldots,z^{i}_{k+n,k+n+1}) \le 0 \text{ along with } z^{i}_{k+n+1,k+n+1} > 0$$

(possibly after implicitly multiplying row k+n+1, $(z^{i}_{1,k+n+1},\ldots,z^{i}_{k+n+1,k+n+1})$, by -1). As already mentioned in discussing "General iteration step" above, a lack of solutions for (Eq) may also be indicated by exhausting $L^{(\mathrm{i})}$ in a MinorP instance or a MajorP instance without resulting in a solution of (Eq).

# 3. Numerical illustrations

Five illustrative examples are presented here. In each example, column indices that are put into $\Pi$ are indicated by asterisks in the CSR table, and a check mark ✓ in $Z^{(i)}$ or $P^{(i)}$ indicates the pivoting column. At the end of each example wherein (Eq) has a solution, that solution is obtained through the guideline stated in "Section 2,2 4 – Stopping step for compact implementation" above.

.

```
Example 1: A simple LP problem
```

*Initialization*

$$\begin{pmatrix} f^T & \\ A & b \end{pmatrix} = \left(\begin{array}{cc|c} -1 & 1 & \\ \hline 1 & 1 & 10 \\ -1 & 0 & -5 \end{array}\right)$$

.

$P^{(0)} =$

| 0 | 0 | 1 | 1 | 10 |
|---|---|---|---|---|
| 0 | 0 | -1 | 0 | -5 |
| -1 | 1 | 0 | 0 | 1 |
| -1 | 0 | 0 | 0 | -1 |
| -10 | 5 | -1 | 1✓ | 0 |

| | Z | P |
|---|---|---|
| 1 | 4 | |

.

*Iteration #1*

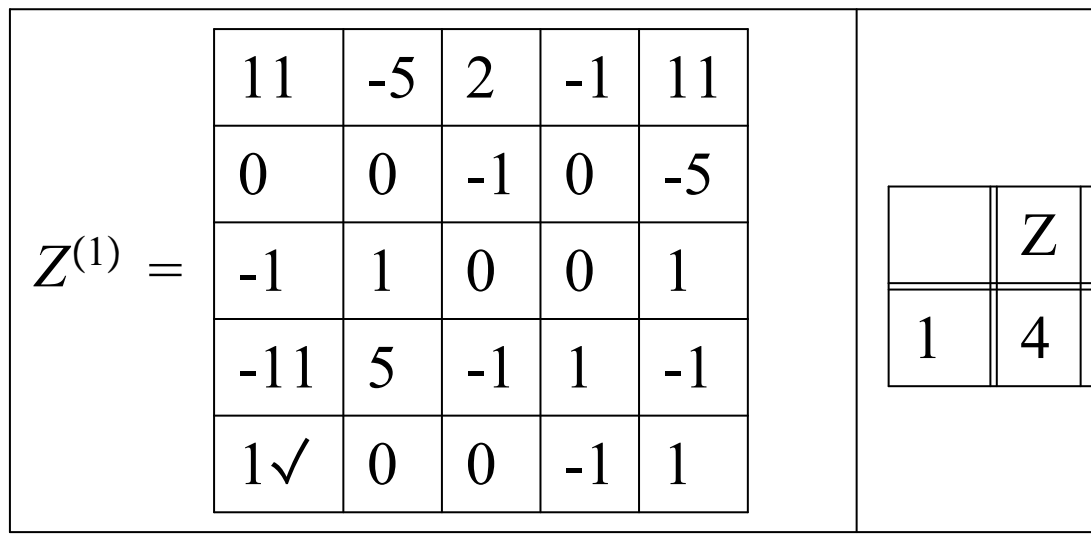

$Z^{(1)} =$

| 11 | -5 | 2 | -1 | 11 |
|---|---|---|---|---|
| 0 | 0 | -1 | 0 | -5 |
| -1 | 1 | 0 | 0 | 1 |
| -11 | 5 | -1 | 1 | -1 |
| 1✓ | 0 | 0 | -1 | 1 |

| | Z | P |
|---|---|---|
| 1 | 4 | 1* |

$P^{(1)} =$

| 0.09 | -0.45 | 0.18 | -0.09 | 1 |
|---|---|---|---|---|
| 0 | 0 | -1 | 0 | -5 |
| 0.09 | 0.55 | 0.18 | -0.09 | 2 |
| 1 | 0 | 1 | 0 | 10 |
| -0.09 | 0.45✓ | -0.18 | -0.91 | 0 |

| | Z | P |
|---|---|---|
| 1 | 4 | 1* |
| 2 | 2* | |

.

*Iteration #2*

$Z^{(2)} =$

| | | | | |
|---|---|---|---|---|
| 0 | 1 | -1 | -1 | -4 |
| -0.2 | 2.2 | -2.6 | -2 | -11 |
| 0.2 | -1.2 | 1.6 | 1 | 8 |
| 1 | 0 | 1 | 0 | 10 |
| 0 | -1 | 1✓ | 0 | 5 |

| | Z | P |
|---|---|---|
| 1 | 4 | 1* |
| 2 | 2* | 3* |

$P^{(2)} =$

| | | | | |
|---|---|---|---|---|
| 0.13 | 0.25 | 0.63 | -0.38 | 1 |
| 0.13 | 0.25 | 1.63 | -0.38 | 2 |
| 0.13 | -0.75 | 0.63 | 0.63 | 5 |
| 0.88 | 0.75 | -0.63 | -0.63 | 5 |
| -013 | -0.25 | -0.63 | -0.63 | 0 |

In the final CSR table, each one of 1,2,3,4 is included an odd number of times. Therefore, a set of solutions of the corresponding LP problem (P) and its dual LP problem is given by

| dual | primal |
|---|---|
| $(1,2)^T$ | $(5,5)^T$ |

.

`Example 2: The instance of Klee-Minty LP problem with n=3`

*Initialization*

$$\begin{pmatrix} f^T & \\ A & b \end{pmatrix} = \left(\begin{array}{ccc|c} 100 & 10 & 1 & \\ \hline 1 & 0 & 0 & 1 \\ 20 & 1 & 0 & 100 \\ 200 & 20 & 1 & 10000 \end{array}\right)$$

.

$P^{(0)} =$

| | | | | | | |
|---|---|---|---|---|---|---|
| 0 | 0 | 0 | 1 | 0 | 0 | 1 |
| 0 | 0 | 0 | 20 | 1 | 0 | 100 |
| 0 | 0 | 0 | 200 | 20 | 1 | 10000 |
| -1 | -20 | -200 | 0 | 0 | 0 | -100 |
| 0 | -1 | -20 | 0 | 0 | 0 | -10 |
| 0 | 0 | -1 | 0 | 0 | 0 | -1 |
| -1 | -100 | -10000 | 100 | 10 | 1✓ | 0 |

| | Z | P |
|---|---|---|
| 1 | 6 | |

.

*Iteration #1*

$Z^{(1)} =$

| 0 | 0 | 0 | 1 | 0 | 0 | 1 |
|---|---|---|---|---|---|---|
| 0 | 0 | 0 | 20 | 1 | 0 | 100 |
| 1 | 100 | 10001 | 100 | 10 | -1 | 10001 |
| -1 | -20 | -200 | 0 | 0 | 0 | -100 |
| 0 | -1 | -20 | 0 | 0 | 0 | -10 |
| -1 | -100 | -10001 | 100 | 10 | 1 | -1 |
| 0 | 0 | 1✓ | 0 | 0 | -1 | 1 |

|  | Z | P |
|---|---|---|
| 1 | 6 | 3* |

.

$P^{(1)} =$ $10^4 *$

| 0 | 0 | 0 | 0.0001 | 0 | 0 | 0.0001 |
|---|---|---|---|---|---|---|
| 0 | 0 | 0 | 0.0020 | 0.0001 | 0 | 0.01 |
| 0.00 | 0.00 | 0.00 | 0.00 | 0.00 | -0.00 | 0.0001 |
| -0.00 | -0.00 | 0.00 | 0.00 | 0.00 | -0.00 | 0.01 |
| 0.00 | -0.00 | 0.00 | 0.00 | 0.00 | -0.00 | 0.0001 |
| 0 | 0 | 0.00 | 0.02 | 0.0002 | 0 | 1 |
| -0.00 | -0.00 | -0.00 | -0.00 | -0.00 | -0.0001 | 0 |

In the final CSR table, each one of 3, 6 is included an odd number of times, whereas each one of 1, 2, 4, 5 is included an even number of times (0 in this case). Therefore, a set of solutions of the corresponding LP problem (P) and its dual LP problem is given by

| dual | primal |
|---|---|
| $(0,0,1)^T$ | $(0,0,10000)^T$ |

It turns out that the Klee-Minty LP problem having *n* variables is also solved in one iteration, regardless of the value of *n*.

.

`Example 3: A very instructive example`

This is an example of a MajorP instance reversing a previous MajorP selection, with the algorithm terminated at that instant in accordance with Claim 7.1 of Section 7.1.1 of [6].

*Initialization*

$$\begin{pmatrix} f^T & \\ A & b \end{pmatrix} = \left(\begin{array}{ccc|c} -9 & 1 & -1 & \\ \hline -2 & -2 & 1 & -7 \\ -4 & 3 & -2 & -3 \end{array}\right)$$

.

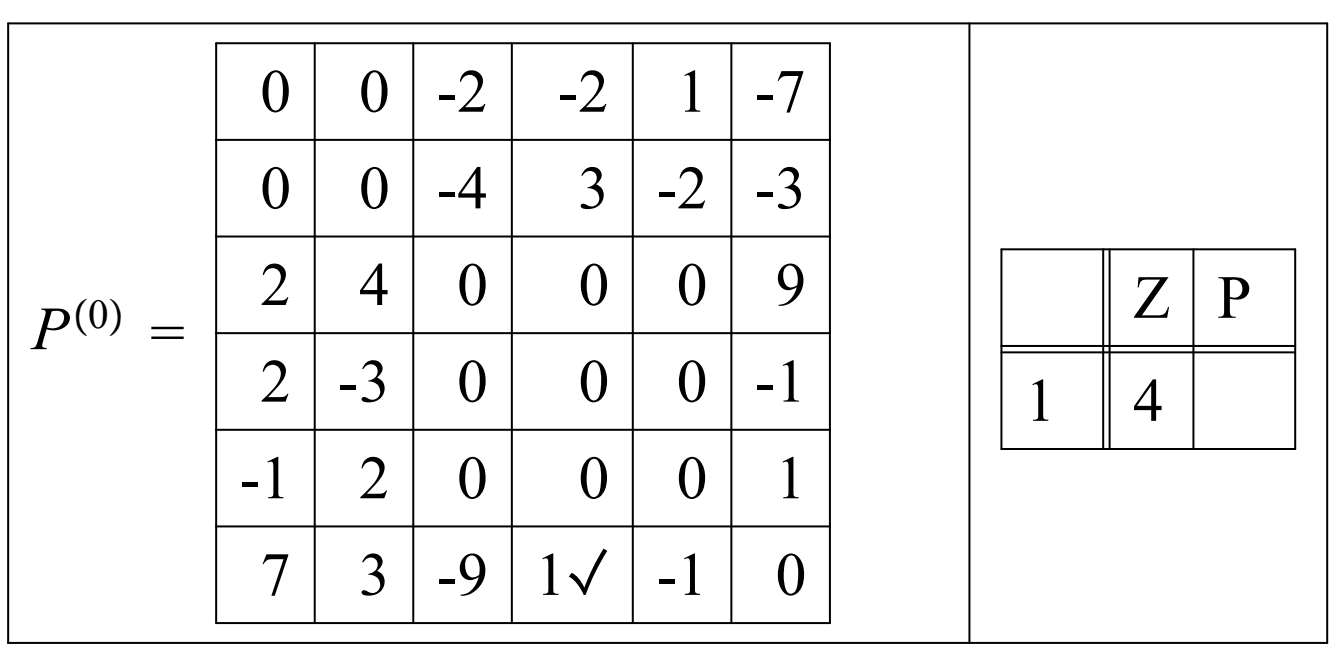

$P^{(0)} =$

| 0 | 0 | -2 | -2 | 1 | -7 |
|---|---|---|---|---|---|
| 0 | 0 | -4 | 3 | -2 | -3 |
| 2 | 4 | 0 | 0 | 0 | 9 |
| 2 | -3 | 0 | 0 | 0 | -1 |
| -1 | 2 | 0 | 0 | 0 | 1 |
| 7 | 3 | -9 | 1√ | -1 | 0 |

| | Z | P |
|---|---|---|
| 1 | 4 | |

.

*Iteration #1*

$Z^{(1)} =$

| 18 | 0 | -20 | 2 | -1 | -9 |
|---|---|---|---|---|---|
| -27 | 0 | 23 | -3 | 1 | 0 |
| 2 | 4 | 0 | 0 | 0 | 9 |
| 9 | 0 | -9 | 1 | -1 | -1 |
| -1 | 2 | 0 | 0 | 0 | 1 |
| -2 | 3√ | 0 | -1 | 0 | 1 |

| | Z | P |
|---|---|---|
| 1 | 4 | 2* |

$P^{(1)} =$

| 18 | 0 | -20 | 2 | -1 | -9 |
|---|---|---|---|---|---|
| -9.67 | 0.33 | 7.67 | -1.33 | 0.33 | 0.33 |
| 40.67 | -1.33 | -30.67 | 5.33 | -1.33 | 7.67 |
| 9 | 0 | -9 | 1 | -1 | -1 |
| 18.33 | -0.67 | -15.33 | 2.67 | -0.67 | 0.33 |
| 27√ | -1 | -23 | 3 | -1 | 0 |

| | Z | P |
|---|---|---|
| 1 | 4 | 2* |
| 2 | 1 | |

.

.

*Iteration #2*

$Z^{(2)} =$

| 0.06 | 0 | -1.11 | 0.11 | -0.06 | -0.5 |
|---|---|---|---|---|---|
| 0.54 | 0.33 | -3.07 | -0.26 | -0.20 | -4.5 |
| -2.26 | -1.33 | 14.52 | 0.81 | 0.93 | 28 |
| -0.50 | 0 | 1 | 0 | -0.5 | 3.5 |
| -1.02 | -0.67 | 5.04 | 0.63 | 0.35 | 9.5 |
| -1.5 | -1 | 7√ | 0 | 0.5 | 13.5 |

| | Z | P |
|---|---|---|
| 1 | 4 | 2* |
| 2 | 1 | 3 |

$P^{(2)} =$

| -0.12 | -0.10 | 0.08 | 0.17 | 0.02 | 1.64 |
|---|---|---|---|---|---|
| 0.06 | 0.05 | 0.21 | -0.09 | -0.01 | 1.43 |
| -0.16 | -0.09 | 0.07 | 0.06 | 0.06 | 1.93 |
| -0.34 | 0.09 | -0.07 | -0.06 | -0.56 | 1.57 |
| -0.23 | -0.20 | -0.35 | 0.35 | 0.03 | -0.21 |
| -0.41 | -0.36 | -0.48 | -0.39 | 0.05✓ | 0 |

| | Z | P |
|---|---|---|
| 1 | 4 | 2 |
| 2 | 1 | 3 |
| 3 | 5* | |

.

*Iteration #3*

$Z^{(3)} =$

| 0.0 | 0.0 | 0.25 | -0.0 | -0.5 | 1.75 |
|---|---|---|---|---|---|
| 0.0 | 0.0 | 0.13 | -0.0 | 0.25 | 1.38 |
| 0.33 | 0.33 | 0.79 | -0.67 | -2.08 | 2.38 |
| -4.67 | -3.67 | -6.46 | 6.33 | 18.42 | -2.38 |
| -7.67 | -6.67 | -11.33 | 11.33 | 32.67 | -7.0 |
| 0.0 | 0.0 | 0.13✓ | -1 | -1.75 | 0.38 |

| | Z | P |
|---|---|---|
| 1 | 4 | 2 |
| 2 | 1 | 3 |
| 3 | 5* | 3* |

$P^{(3)} =$

| -0.11 | -0.11 | -0.32 | 0.21 | 0.16 | 1 |
|---|---|---|---|---|---|
| -0.05 | -0.05 | -0.16 | 0.11 | 0.58 | 1 |
| 0.42 | 0.42 | 1.26 | -0.84 | -2.63 | 3 |
| -1.95 | -0.95 | 8.16 | 0.89 | 1.42 | 17 |
| -2.89 | -1.89 | 14.32 | 1.79 | 2.84 | 27 |
| -0.05 | -0.05 | -0.16 | -0.89 | -1.42 | 0 |

In the final CSR table, each one of 1,2,4,5 is included an odd number of times, whereas the number 3 is included an even number of times. Therefore, a set of solutions of the corresponding LP problem (P) and its dual LP problem is given by

| dual | primal |
|---|---|
| $(1,1)^T$ | $(0,17,27)^T$ |

.

`Example 4: An infeasible LP problem`

*Initialization*

$$\begin{pmatrix} f^T & \\ A & b \end{pmatrix} = \left(\begin{array}{cc|c} 1 & 1 & \\ \hline -1 & 2 & -4 \\ 2 & 1 & 3 \end{array}\right)$$

.

$P^{(0)} =$

| | | | | |
|---|---|---|---|---|
| 0 | 0 | -1 | 2 | -4 |
| 0 | 0 | 2 | 1 | 3 |
| 1 | -2 | 0 | 0 | -1 |
| -2 | -1 | 0 | 0 | -1 |
| 4 | -3 | 1√ | 1 | 0 |

| | Z | P |
|---|---|---|
| 1 | 3 | |

.

*Iteration #1*

$Z^{(1)} =$

| | | | | |
|---|---|---|---|---|
| 5 | -5 | 1 | 3 | -5 |
| -10 | 10 | -2 | -1 | 5 |
| 5 | -5 | 1 | 1 | -1 |
| -2 | -1 | 0 | 0 | -1 |
| -1 | 2√ | -1 | 0 | 1 |

| | Z | P |
|---|---|---|
| 1 | 3 | 2* |

$P^{(1)} =$

| | | | | |
|---|---|---|---|---|
| 0 | 0.5 | 0 | 2.5 | -2.5 |
| -1 | 0.1 | -0.2 | -0.1 | 0.5 |
| 0 | 0.5 | 0 | 0.5 | 1.5 |
| -3 | 0.1 | -0.2 | -0.1 | -0.5 |
| 1 | -0.2 | -0.6 | 0.2√ | 0 |

| | Z | P |
|---|---|---|
| 1 | 3 | 2* |
| 2 | 4 | |

.

*Iteration #2*

$Z^{(2)} =$

| | | | | |
|---|---|---|---|---|
| -75 | 3 | -5 | 25 | -15 |
| 2 | 0 | 0 | -1 | 1 |
| -15 | 1 | -1 | 5 | -1 |
| 30 | -1 | 2 | -10 | 5 |
| -5√ | 0 | -1 | 2 | -1 |

| | Z | P |
|---|---|---|
| 1 | 3 | 2* |
| 2 | 4 | 1 |

$P^{(2)} =$

| | | | | |
|---|---|---|---|---|
| -0.01 | -0.04 | 0.07 | -0.33 | 0.2 |
| 0.03 | 0.08 | -0.13 | -0.33 | 0.6 |
| -0.2 | 0.4 | 0 | 0 | 2.0 |
| 0.40 | 0.20 | 0 | 0 | -1.0 |
| -0.07 | -0.20 | -0.67 | 0.33√ | 0 |

| | Z | P |
|---|---|---|
| 1 | 3 | 2 |
| 2 | 4 | 1 |
| 3 | 4* | |

.

*Iteration #3*

$Z^{(3)} =$

| 0.32 | -0.04 | -0.60 | 1 | -0.8 |
|---|---|---|---|---|
| 0.36 | 0.08 | -0.80 | 1 | -0.4 |
| -0.20 | 0.40 | 0 | 0 | 2 |
| 1 | -0.00 | -2 | 3 | -3 |
| -0.40 | -0.20 | 0 | -1 | 1 |

| | Z | P |
|---|---|---|
| 1 | 3 | 2 |
| 2 | 4 | 1 |
| 3 | 4* | na |

From the last row of $Z^{(3)}$, it is clear that corresponding (Eq) has no solutions. Note that, for this particular example, one could choose columns along the way to have the CSR

| | Z | P |
|---|---|---|
| 1 | 4 | 1 |
| 2 | 3 | 2 |
| 3 | 4 | na |

instead of

| | Z | P |
|---|---|---|
| 1 | 3 | 2 |
| 2 | 4 | 1 |
| 3 | 4 | na |

.

`Example 5: Another very instructive example`

This is an example of a MinorP instance having to reverse a previous MajorP selection (in iteration 5).

*Initialization*

$$\begin{pmatrix} f^T & \\ A & b \end{pmatrix} = \left(\begin{array}{cccc|c} 3 & 4 & 1 & 7 & \\ \hline 8 & 3 & 4 & 1 & 7 \\ 2 & 6 & 1 & 5 & 3 \\ 1 & 4 & 5 & 2 & 8 \end{array}\right)$$

.

$P^{(0)} =$

| 0 | 0 | 0 | 8 | 3 | 4 | 1 | 7 |
|---|---|---|---|---|---|---|---|
| 0 | 0 | 0 | 2 | 6 | 1 | 5 | 3 |
| 0 | 0 | 0 | 1 | 4 | 5 | 2 | 8 |
| -8 | -2 | -1 | 0 | 0 | 0 | 0 | -3 |
| -3 | -6 | -4 | 0 | 0 | 0 | 0 | -4 |
| -4 | -1 | -5 | 0 | 0 | 0 | 0 | -1 |
| -1 | -5 | -2 | 0 | 0 | 0 | 0 | -7 |
| -7 | -3 | -8 | 3 | 4 | 1✓ | 7 | 0 |

| | Z | P |
|---|---|---|
| 1 | 6 | |

The algorithm terminated in five iterations with a solution of (Eq). The final CSR in this case is

.

| | Z | P |
|---|---|---|
| 1 | 6 | 3 |
| 2 | 4 | 1 |
| 3 | 5 | 2 |
| 4 | 7 | 5* |
| 5 | 3* | 6* |

In the final CSR table, each one of 1,2,4,7 is included an odd number of times, whereas each one of 3,5,6 is included an even number of times. Therefore, a set of solutions of the corresponding LP problem (P) and its dual LP problem is given by

| dual | primal |
|---|---|
| $(0.0263, 1.3947, 0)^T$ | $(0.8421, 0, 0, 0.2632)^T$ |

# 4. Some directions for further work

The rule in the underlying algorithm [6] by which MinorP pivoting makes its column selections is translated in this compact implementation article into the rule by which elements of $L^{(i)}$ are ordered (recall the definition $L^{(i)} = \{j \text{ such that } p^{i-1}_{k+n+1,j} > 0\}$). That rule has some freedom about it, while still maintaining required agreement with stated computational complexity in Section 7 of [6].

For example, the computational complexity statement allows one to change the rule by which elements of $L^{(i)}$ are ordered to "ascending order of $j$ having $p^{i-1}_{k+n+1,j} > 0$" (that is, leftmost $j$), in place of "ascending order of $p^{i-1}_{k+n+1,j} > 0$" (least positive entry in the last row). We find that as a consequence of that change, *Example 5* LP problem (above) is solved in 3 iterations instead of the 5 iterations reported above, with corresponding CSR being

.

| | Z | P |
|---|---|---|
| 1 | 4 | 1 |
| 2 | 5 | 2 |
| 3 | 7* | 5* |

in place of

| | Z | P |
|---|---|---|
| 1 | 6 | 3 |
| 2 | 4 | 1 |
| 3 | 5 | 2 |
| 4 | 7 | 5* |
| 5 | 3* | 6* |

But that new rule does not solve *Example 2* (instance of Klee-Minty LP problem with n=3) LP problem in 1 iteration; it solves *Example 2* in 2 iterations.

Thus, one may surmise that different rules for ordering elements of $L^{(i)}$ in MinorP instances are well-suited to different classes of practical problems, while still maintaining the algorithm's strongly polynomial-time computational complexity. Accordingly, one direction for further work related to this article is to perform extensive numerical computations to discover different efficacious rules for different classes of practical LP problems.

Another direction for further work is to investigate connections between the complementary pivoting presented in this article and the general complementary pivoting presented in the article [4], with "Pivot-in" and "Pivot-out" corresponding to MinorP and MajorP respectively. This kind of investigation may yield new useful rules for ordering elements of the set $L^{(i)}$.

This compact implementation generates a good deal of numerical data during intermediate iterations in solving each LP problem. As a direction for further work, one may want to investigate how useful the intermediate data (generated by this compact implementation) can be when the algorithm is used as a subroutine for solving certain integer programming problems.